\documentclass[12pt,a4paper]{amsart}
\usepackage{mathrsfs}

\newtheorem{theo+}              {Theorem}           [section]
\newtheorem{prop+}  [theo+]     {Proposition}
\newtheorem{coro+}  [theo+]     {Corollary}
\newtheorem{lemm+}  [theo+]     {Lemma}
\newtheorem{exam+}  [theo+]     {Example}
\newtheorem{rema+}  [theo+]     {Remark}
\newtheorem{defi+}  [theo+]     {Definition}

\newenvironment{theorem}{\begin{theo+}}{\end{theo+}}

\newenvironment{corollary}{\begin{coro+}}{\end{coro+}}
\newenvironment{lemma}{\begin{lemm+}}{\end{lemm+}}

\usepackage{amsthm}
\theoremstyle{plain} \theoremstyle{remark}

\newtheorem{example}{Example}

\newtheorem*{ack}{\bf Acknowledgments}

\def \r{\mbox{${\mathbb R}$}}
\def\E{/\kern-1.0em \equiv }

\title{On the generalized Chen's conjecture on biharmonic submanifolds }
\author{Ye-Lin Ou$^{*}$ and Liang Tang}
\address{Department of
Mathematics,\newline\indent Texas A $\&$ M University-Commerce,
\newline\indent Commerce TX 75429,\newline\indent USA.\newline\indent
E-mail:yelin$\_$ou@tamu-commerce.edu \newline\indent \vskip0.1cm
School of Mathematics and Computer Science,\newline\indent Guangxi
University for Nationalities,
\newline\indent 188 East Daxue Road,\newline\indent Nanning, Guangxi 530006,\newline\indent P. R. China}
\thanks{ \indent* Supported by Texas A $\&$ M University-Commerce ``Faculty
Research Enhancement Project" (2010-11)}

\begin{document}
\title[Generalized Chen's conjecture on biharmonic submanifolds is false]{On the generalized Chen's conjecture on biharmonic submanifolds}
\date {01/18/2011} 
 \maketitle

\section*{Abstract}
\begin{quote}
{\footnotesize The generalized Chen's conjecture on biharmonic
submanifolds asserts that any biharmonic submanifold of a
non-positively curved manifold is minimal. In this paper, we prove
that this conjecture is false by constructing a foliation of proper
biharmonic hyperplanes in a $5$-dimensional conformally flat space
with negative sectional curvature. Many examples of proper
biharmonic submanifolds of non-positively curved spaces are also
given.}
\end{quote}
\section{Biharmonic submanifolds and the generalized Chen's conjecture}

We work on the category of smooth objects so all manifolds, maps, and tensor fields discussed in this paper are smooth unless there is an otherwise statement.\\

A {\bf biharmonic map} is a map $\varphi:(M, g)\longrightarrow (N,
h)$ between Riemannian manifolds which is a solution of the $4$th
order PDEs
\begin{equation}\label{BTF}
\tau^{2}(\varphi):={\rm
Trace}_{g}(\nabla^{\varphi}\nabla^{\varphi}-\nabla^{\varphi}_{\nabla^{M}})\tau(\varphi)
- {\rm Trace}_{g} R^{N}({\rm d}\varphi, \tau(\varphi)){\rm d}\varphi
=0,
\end{equation}
 where $R^{N}$ denotes the curvature operator of $(N, h)$ defined by
$$R^{N}(X,Y)Z=
[\nabla^{N}_{X},\nabla^{N}_{Y}]Z-\nabla^{N}_{[X,Y]}Z,$$ and
 $\tau(\varphi)={\rm Trace}_{g}\nabla {\rm d} \varphi$ is the tension
field of $\varphi$ and $\tau(\varphi)=0$ means the map $\varphi$ is
harmonic.\\

Clearly, it follows from (\ref{BTF}) that any harmonic map is
biharmonic and we
call those non-harmonic biharmonic maps {\bf proper biharmonic maps}.\\

A submanifold $M$ of $(N,h)$ is called a {\bf biharmonic
submanifold} if the inclusion map ${\bf i}: (M, {\bf
i}^*h)\longrightarrow (N,h)$ is a biharmonic isometric immersion. It
is well-known that an isometric immersion is minimal if and only if
it is harmonic. So a minimal submanifold is trivially biharmonic and
we call a non-minimal biharmonic submanifold a {\bf proper
biharmonic submanifold}.\\

Among the problems in the study of biharmonic maps the following are
two fundamental ones: {\bf (1) existence problem:} given two model
spaces (e.g., some ``good" spaces such as spaces of constant
sectional curvature or more general symmetric or homogeneous
spaces), does there exist a proper biharmonic map from one space
into another? {\bf (2) classification problem:} classify all proper
biharmonic maps between two model spaces where the existence is
known. A typical and
challenging classification problem is the following \\

{\bf Chen's conjecture \cite{CH}:} any biharmonic submanifold in
Euclidean
space is minimal.\\

The conjecture was proved to be true for biharmonic surfaces in
$\r^{3}$ (by Jiang \cite{Ji2} and Chen-Ishikawa \cite{CI}
independently) and for biharmonic hypersurface in $\r^4$
(\cite{HV}). Dimitri$\acute{\rm c}$ \cite{Di} showed that the
conjecture is also true for any biharmonic curves, any biharmonic
submanifold of finite type, any pseudo-umbilical submanifolds
$M^m\subset \r^n$ with $m\ne 4$, and  any biharmonic hypersurface in
$\r^n$ with at most two distinct principal curvatures. However, the
conjecture is still open in general.\\

In the same direction of classifying proper biharmonic submanifolds
of non-positively curved manifolds, Caddeo, Montaldo and Oniciuc
\cite{CMO2} shown that any biharmonic submanifold in hyperbolic
$3$-space $H^{3}(-1)$ is minimal, and any pseudo-umbilical
biharmonic submanifold $M^m\subset H^n$ with $m\ne 4$ is minimal. It
is also shown in \cite{BMO1} that any biharmonic hypersurface of
$H^n$ with at most two distinct principal curvatures is minimal. All
these results suggest the following generalized Chen's conjecture on
biharmonic submanifolds which was proposed by Caddeo, Montaldo and
Oniciuc \cite{CMO1}:\\

{\bf The generalized Chen's conjecture}: any biharmonic submanifold
of $(N, h)$ with ${\rm Riem}^N\leq 0$ is minimal (see e.g.,
\cite{CMO1}, \cite{MO}, \cite{BMO1}, \cite{BMO2}, \cite{BMO3},
\cite{Ba1}, \cite{Ba2}, \cite{Ou1}, \cite{Ou2}, \cite{IIU}).\\

The gaol of this paper is to  prove that the generalized Chen's
conjecture for biharmonic submanifolds is false. We accomplished
this by using the idea of constructing foliations of proper
biharmonic hyperplanes in a conformally flat space given in
\cite{Ou1}. The idea is to determine a conformally flat metric on
$\r^{m+1}$ so that a foliation by the hyperplanes defined by the
graphs of a linear function becomes a proper biharmonic foliation.
It turns out that when $m=4$ the system of biharmonicity equations
reduces to a single equation which has many solutions including
counter examples to the generalized Chen's conjecture.

\section{Foliations of conformally flat spaces by biharmonic hyperplanes}

\indent Two Riemannian metrics $g$ and $\bar{g}$ on $M$ are said to
be conformally equivalent, if there exists a function $\sigma$ on
$M$ such that $\bar{g}= e^{2\sigma}g$. A map $\varphi:(M,
g)\longrightarrow (N, h)$ between Riemannian manifolds is said to be
conformal if there exists a function $\sigma$ on $M$ such that
${\varphi}^{*} h= e^{2\sigma}g$. Two Riemannian manifolds ($M,g$)
and ($N,h$) are said to be conformally diffeomorphic, if there
exists a conformal diffeomorphism $\varphi:(M, g)\longrightarrow (N,
h)$. A Riemannian manifold ($M^{m},g$) is called a conformally flat
space if for any point of $M$ there exists a neighborhood which is
conformally diffeomorphic to an open neighborhood of the Euclidean
space $\mathbf{R}^{m}$. Notice that any two-dimensional Riemannian
manifold is conformally flat due to the existence of isothermal
coordinates. For $m= 3$, ($M^{m},g$) is conformally flat if and only
if the Schouten tensor $H$ satisfies $({\nabla}_{X} H)(Y,Z) =
({\nabla}_{Y}H)(X,Z)$ for any vector fields $X,\;Y$ and $Z$ on $M$.
For $m\geq 4$, $(M^{m},g)$ is conformally flat if and only if the
Weyl curvature vanishes identically. It is well known that a space
of constant sectional curvature is conformally flat but there exist
many conformally flat
spaces which are not of constant sectional curvature.\\
\indent Let $\nabla, \rm R, \rm Ric, r, \rm K$ (respectively ${\bar \nabla},
\bar {\rm R}, {\overline {\rm Ric}}, {\bar r}, \bar {\rm K}$) denote the Levi-Civita
connection, Riemannian curvature, Ricci curvature, scalar curvature,
and sectional curvature of the Riemannian metric $g$ (respectively
${\bar g}=e^{2\sigma}g$). Then, we have
\begin{eqnarray}\notag
{\bar R}(W,Z, X,Y)&=&e^{2\sigma}\{R(W,Z,
X,Y)+g(\nabla_{X}\nabla\sigma,Z)g(Y,W)\\\label{c}
&&-g(\nabla_{Y}\nabla\sigma,Z)g(X,W)+g(X,Z)g(\nabla_{Y}\nabla
\sigma,W)\\\notag &&-g(Y,Z)g(\nabla_{X}\nabla \sigma,W)
+[(Y\sigma)(Z\sigma)-g(Y,Z)|\nabla \sigma|^{2}]g(X,W)\\\notag &&
-[(X\sigma)(Z\sigma)-g(X,Z)|\nabla \sigma|^{2}]g(Y,W)\\\notag
&&+[(X\sigma)g(Y,Z)-(Y\sigma)g(X,Z)]g(\nabla \sigma,W)\}.
\end{eqnarray}

In local coordinates, (replacing $W, Z, X, Y$ by
$\partial_i,
\partial_j, \partial_k, \partial_l,$ respectively in (\ref{c}) to get)
\begin{eqnarray}\notag\label{d}
e^{-2\sigma}{\bar
R}_{ij\,kl}&=&R_{ij\,kl}+g_{il}\sigma_{jk}-g_{ik}\sigma_{jl}+g_{jk}
\sigma_{il}-g_{jl}\sigma_{ik}\\\label{246}&&
+(g_{il}g_{jk}-g_{ik}g_{jl})|\nabla \sigma|^{2},
\end{eqnarray}
where we have used the  notation
$\sigma_{jl}=\nabla_{l}\sigma_j-\sigma_l\sigma_j=\nabla_{l}\nabla_j\sigma-\sigma_l\sigma_j$.\\
and
\begin{eqnarray}\label{e}
{\bar R}_{jk}=R_{jk}-(n-2)\sigma_{jk}-g_{jk}[ \Delta\sigma
+(n-2)|\nabla \sigma|^{2}].
\end{eqnarray}

\indent Let $P$ be a section spanned by an orthonormal basis $X, Y$
with respect to $g$ (${\bar X}=e^{-\sigma}X, {\bar Y}=e^{-\sigma}Y$
form an orthonormal basis w.r.t. ${\bar g}$). Then, we have the
relationship between the sectional curvatures with respect to
metrics ${\bar g}$ and $g$ given by
\begin{eqnarray}\label{QL}
e^{2\sigma}{\bar K}(P)=K(P)-(g(\nabla_{X}\nabla
\sigma,X)+g(\nabla_{Y}\nabla \sigma,Y))\\\notag
-(|\nabla\sigma|^{2}-(X\sigma )^{2}-(Y\sigma)^{2}).
\end{eqnarray}

\begin{theorem}\label{EQ1}\cite{Ou1}
Let $\varphi:M^{m}\longrightarrow N^{m+1}$ be an isometric immersion
of codimension-one with mean curvature vector $\eta=H\xi$. Then
$\varphi$ is biharmonic if and only if:
\begin{equation}\label{eq1}
\begin{cases}
\Delta_g H-H |A|^{2}+H{\rm
Ric}^N(\xi,\xi)=0,\\
 2A\,({\rm grad}_g\,H) +\frac{m}{2} {\rm grad}_g\, H^2
-2\, H \,({\rm Ric}^N\,(\xi))^{\top}=0,
\end{cases}
\end{equation}
where  $\Delta_g$ and ${\rm grad}_g$ are the Laplacian and the
gradient operators of the hypersurface, and ${\rm Ric}^N :
T_qN\longrightarrow T_qN$ denotes the Ricci operator of the ambient
space defined by $\langle {\rm Ric}^N\, (Z), W\rangle={\rm Ric}^N
(Z, W)$ and  $A$ is the shape operator of the hypersurface with
respect to the unit normal vector $\xi$.
\end{theorem}

A further computation yields
\begin{eqnarray}\label{j}
\rm Ric(\xi, \xi)&=& e^{2\sigma}c_{m+1}^{j}c_{m+1}^{k}\rm Ric_{jk}\\\notag
&=&e^{2\sigma}c_{m+1}^{j}c_{m+1}^{k}\{(m-1)\rm
Hess(\sigma)(\partial_j,
\partial_k)\\\notag &&-(m-1)\sigma_k\sigma_j +e^{-2\sigma}\delta_{jk}[
\Delta\sigma +(m-1)|\rm grad\sigma|^{2}]\}\\\notag
&=&\Delta\sigma+(m-1)[\rm Hess(\sigma)(\xi, \xi)-\xi^2(\sigma)+|\rm
grad\sigma|^{2}].
\end{eqnarray}
\begin{eqnarray}\label{k}
[\rm Ric(\xi)]^T&=&
\sum_{i=1}^{m}e^{2\sigma}c_{m+1}^{j}c_{i}^{k}\rm Ric_{jk}e_i\\\notag
&=&\sum_{i=1}^{m}e^{2\sigma}c_{m+1}^{j}c_{i}^{k}\{(m-1)\rm
Hess(\sigma)(\partial_j,
\partial_k)\\\notag &&-(m-1)\sigma_k\sigma_j +e^{-2\sigma}\delta_{jk}[
\Delta\sigma +(m-1)|\rm grad\sigma|^{2}]\}e_i\\\notag
&=&\sum_{i=1}^{m}(m-1)[\rm Hess(\sigma)(\xi,
e_i)e_i-\xi(\sigma)e_i(\sigma)e_i]\\\notag
&=&(m-1)[\sum_{i=1}^{m}\rm Hess(\sigma)(\xi, e_i)e_i-\xi(\sigma){\rm
grad}_g \sigma]
\\\notag &=&(m-1)[\sum_{i=1}^{m}\{e_i(\xi(\sigma))e_i-(\nabla_{e_i}\xi)(\sigma)e_i\}-\xi(\sigma){\rm
grad}_g
\sigma]\\\notag &=&(m-1)[{\rm
grad}_g\xi(\sigma)+\sum_{i=1}^{m}
(Ae_i)(\sigma)e_i-\xi(\sigma){\rm
grad}_g \sigma]
\\\notag &=&(m-1)[{\rm
grad}_g\xi(\sigma)-\xi(\sigma){\rm
grad}_g\sigma+\sum_{i,j=1}^{m}b(e_i,e_j)e_j(\sigma)e_i]
\\\notag &=&(m-1)[{\rm
grad}_g\xi(\sigma)-\xi(\sigma){\rm
grad}_g\sigma+\sum_{i,j=1}^{m}b(e_i,e_j(\sigma)e_j)e_i]
\\\notag &=&(m-1)[{\rm
grad}_g\xi(\sigma)-\xi(\sigma){\rm
grad}_g\sigma+\sum_{i=1}^{m}b(e_i,{\rm
grad}_g\sigma)e_i]
\\\notag &=&(m-1)[{\rm
grad}_g\xi(\sigma)-\xi(\sigma){\rm
grad}_g\sigma+\sum_{i=1}^{m}h (A({\rm
grad}_g
\sigma),e_i)e_i]
\\\notag &=&(m-1)[{\rm
grad}_g\xi(\sigma)-\xi(\sigma){\rm
grad}_g\sigma+A({\rm
grad}_g \sigma)].
\end{eqnarray}

Using the following conventions for indices:\\
$i, j, k,l=1, 2, \cdots, m;\;\alpha,\beta=1, 2,\cdots, m, m+1$.\\

\begin{theorem}\label{f1}
Let $a_i, i=1, 2, \ldots, m$ and $c$ be constants. Then, the
isometric immersion $\varphi : \mathbb{R}^m\longrightarrow
(\mathbb{R}^{m+1},h=f^{-2}(z)(\sum_{i=1}^m{\rm d}{x_i}^{2}+{\rm
d}{z}^{2})$ with $\varphi(x_1,\ldots, x_m)=(x_1,\ldots,
x_m,\sum_{i=1}^{m}a_ix_i+c)$ is biharmonic  if and only if one of
the following three cases happens
\begin{itemize}
\item[(1)] $f'=0$,
in this case $\varphi$ is minimal (actually, totally geodesic), or
\item[(2)]  $m=4$ and
$f$ is a solution of the equation
\begin{equation}\label{single}
\sum_{i=1}^4a_i^2f^2f'''+(4-\sum_{i=1}^4a_i^2)ff'f''-4(2+\sum_{i=1}^4a_i^2)(f')^3=0,
\end{equation}
or
\item[(3)] $a_i=0,\;i=1,\;\cdots,\;m \;{\rm and}\;
f(z)=\frac{1}{Az+B}$, where  $A$ and $B$ are constants. In this case
each hyperplane is a proper biharmonic hypersurface. This recovers a
result (Theorem 3.1) obtained earlier in \cite{Ou1}.
\end{itemize}
\end{theorem}

\begin{proof}
We have $\sigma=\ln f(z)$.Using the notations
$\partial_i=\frac{\partial}{\partial x_i},i=1,2, \ldots, m,
 \partial_{m+1}=\frac{\partial}{\partial z},$ \;we can easily check that the
induced metric is given by
\begin{equation}\notag
g_{ij}=g(\partial_i,\partial_j)=h(d\varphi(\partial
_i),d\varphi(\partial _j))\circ\varphi=
\begin{cases}
(1+a_i^2)f^{-2}(\sum_{i=1}^{m}a_ix_i+c),\;\;i=j,\\\notag
a_ia_jf^{-2}(\sum_{i=1}^{m}a_ix_i+c)\;\;\;\;\;i\ne j.
\end{cases}
\end{equation}
It is easy to check that $\bar e_\alpha=f(z)\partial_\alpha\;$
constitute an orthonormal frame  on $\mathbb{R}^{m+1}$, then we can also check that $e_i=-a_ik_ik_{i-1}\sum_{l=1}^{i-1}a_l\bar
e_l+\frac{k_i}{k_{i-1}}\bar e_i+a_ik_ik_{i-1}\bar e_{m+1}\;\;
(i=1,2,\ldots, m),e_{m+1}=\sum_{l=1}^{m}a_lk_{m}\bar e_l-k_{m}\bar
e_{m+1}$ is an adapted frame tangent to the hypersurface $z=\sum_{i=1}^{m}a_ix_i+c$ with
$\xi=e_{m+1}$ being the unit normal vector field,where
$k_i=\frac{1}{\sqrt{1+\sum_{l=1}^{i}a_l^2}},\;i=1, \cdots,\;m,\;k_0=1,\;a_0=0$. In
fact, we have
\begin{eqnarray*}\notag
\langle e_i,e_i\rangle&=&\langle
-a_ik_ik_{i-1}\sum_{l=1}^{i-1}a_l\bar e_l+\frac{k_i}{k_{i-1}}\bar
e_i+a_ik_ik_{i-1}\bar e_{m+1},\\\notag
&&-a_ik_ik_{i-1}\sum_{l=1}^{i-1}a_l\bar e_l+\frac{k_i}{k_{i-1}}\bar
e_i+a_ik_ik_{i-1}\bar e_{m+1}\rangle\\\notag
&=&a_i^2k_i^2k_{i-1}^2\sum_{l=1}^{i-1}a_l^2+\frac{k_i^2}{k_{i-1}^2}+a_i^2k_i^2k_{i-1}^2\\\notag
&=&a_i^2k_i^2k_{i-1}^2(\sum_{l=1}^{i-1}a_l^2+1)+\frac{k_i^2}{k_{i-1}^2}=a_i^2k_i^2+\frac{k_i^2}{k_{i-1}^2}\\\notag
&=&k_i^2(a_i^2+\frac{1}{k_{i-1}^2})=k_i^2\frac{1}{k_i^2}=1,\\\notag
\langle e_{m+1},e_{m+1}\rangle&=&\langle \sum_{l=1}^{m}a_lk_{m}\bar
e_l-k_{m}\bar e_{m+1},\sum_{l=1}^{m}a_lk_{m}\bar e_l-k_{m}\bar
e_{m+1} \rangle\\\notag
&=&\sum_{l=1}^{m}a_l^2k_{m}^2+k_{m}^2=(\sum_{l=1}^{m}a_l^2+1)k_{m}^2=1,
\end{eqnarray*}
because $\langle e_i,e_j\rangle=\langle e_j,e_i\rangle$,let $i<j$, we have
\begin{eqnarray*}\notag
\langle e_i,e_j\rangle&=&\langle
-a_ik_ik_{i-1}\sum_{l=1}^{i-1}a_l\bar e_l+\frac{k_i}{k_{i-1}}\bar
e_i+a_ik_ik_{i-1}\bar e_{m+1},\\\notag
&&-a_jk_jk_{j-1}\sum_{l=1}^{j-1}a_l\bar e_l+\frac{k_j}{k_{j-1}}\bar
e_j+a_jk_jk_{j-1}\bar e_{m+1}\rangle\\\notag
&=&a_ia_jk_ik_jk_{i-1}k_{j-1}\sum_{l=1}^{i-1}a_l^2-\frac{a_ia_jk_ik_jk_{j-1}}{k_{i-1}}+a_ia_jk_ik_jk_{i-1}k_{j-1}\\\notag
&=&a_ia_jk_ik_jk_{i-1}k_{j-1}(\sum_{l=1}^{i-1}a_l^2+1)-\frac{a_ia_jk_ik_jk_{j-1}}{k_{i-1}}\\\notag
&=&\frac{a_ia_jk_ik_jk_{i-1}k_{j-1}}{k_{i-1}^2}-\frac{a_ia_jk_ik_jk_{j-1}}{k_{i-1}}=0,\\\notag
\langle e_i,e_{m+1}\rangle&=&\langle
-a_ik_ik_{i-1}\sum_{l=1}^{i-1}a_l\bar e_l+\frac{k_i}{k_{i-1}}\bar
e_i+a_ik_ik_{i-1}\bar e_{m+1},\\\notag &&\sum_{l=1}^{m}a_lk_{m}\bar
e_l-k_{m}\bar e_{m+1}\rangle\\\notag
&=&-a_ik_ik_{i-1}k_{m}\sum_{l=1}^{i-1}a_l^2+\frac{a_ik_ik_{m}}{k_{i-1}}-a_ik_ik_{i-1}k_{m}\\\notag
&=&-a_ik_ik_{i-1}k_{m}(\sum_{l=1}^{i-1}a_l^2+1)+\frac{a_ik_ik_{m}}{k_{i-1}}\\\notag
&=&-a_ik_ik_{i-1}k_{m}\frac{1}{k_{i-1}^2}+\frac{a_ik_ik_{m}}{k_{i-1}}=0.
\end{eqnarray*}
\indent it is well-known that $(\mathbb{R}^{m+1},\bar
h=\sum_{i=1}^m{\rm d}{x_i}^{2}+{\rm d}{z}^{2})$ is a flat space, so
$\bar \nabla_{\partial_\alpha}{\partial_\beta}=0$,where $\bar
\nabla$ denotes the Levi-Civita connection with respect to $\bar h$.
using
\begin{eqnarray*}
{\bar \nabla}_{X}Y=\nabla_{X}Y+(X\sigma)Y+(Y\sigma)X-h(X,Y){\rm
grad} \sigma
\end{eqnarray*}
we have the coefficients of the Levi-Civita connection with respect
to $h$ given by
\begin{eqnarray*}\notag
\nabla_{\partial_\alpha}{\partial_\beta}&=&h(\partial_\alpha,\partial_\beta){\rm
grad}
\sigma-\partial_\alpha(\sigma)\partial_\beta-\partial_\beta(\sigma)\partial_\alpha,\\\notag
\nabla_{\bar e_\alpha}{\bar
e_\beta}&=&e^{2\sigma}\nabla_{\partial_\alpha}{\partial_\beta}+e^{\sigma}\partial_\alpha(e^{\sigma})\partial_\beta
=\delta_{\alpha\beta}{\rm
grad} \sigma-\bar e_\beta(\sigma)\bar
e_\alpha,\\\notag {\rm
grad}\sigma&=&\bar e_{m+1}(\sigma)\bar
e_{m+1}=f^{'}\bar e_{m+1},
\end{eqnarray*}
so we have
\begin{equation}
(\nabla_{\bar e_\alpha}\bar e_{\beta})=\left(\begin{array}{ccccc}
f'\bar e_{m+1} &0&\ldots&0&-f'\bar e_1
\\ 0 & f'\bar e_{m+1}&\ldots&0&-f'\bar e_2\\
\ldots&\ldots&\ldots&\ldots&\ldots
\\ 0 & 0&\ldots&f'\bar e_{m+1}&-f'\bar e_m\\
0&0&\ldots&0&0\\
\end{array}\right)_{(m+1)\times (m+1)},
\end{equation}
and
\begin{eqnarray*}\notag
\nabla_{e_i}e_i&=&k_i^2k_{i-1}^2\nabla_{(-a_i\sum_{l=1}^{i-1}a_l\bar
e_l+\frac{1}{k_{i-1}^2}\bar e_i+a_i\bar
e_{m+1})}(-a_i\sum_{l=1}^{i-1}a_l\bar e_l+\frac{1}{k_{i-1}^2}\bar
e_i+a_i\bar e_{m+1})\\\notag
&=&k_i^2k_{i-1}^2\{a_i^2\sum_{l=1}^{i-1}a_l^2\nabla_{\bar e_l}\bar
e_l-a_i^2\sum_{l=1}^{i-1}a_l\nabla_{\bar e_l}\bar
e_{m+1}+\frac{1}{k_{i-1}^4}\nabla_{\bar e_i}\bar
e_i+\frac{a_i}{k_{i-1}^2}\nabla_{\bar e_i}\bar e_{m+1}\}\\\notag
&=&k_i^2k_{i-1}^2f'[(a_i^2\sum_{l=1}^{i-1}a_l^2+\frac{1}{k_{i-1}^4})\bar
e_{m+1}+a_i^2\sum_{l=1}^{i-1}a_l\bar e_l-\frac{a_i}{k_{i-1}^2}\bar
e_i],
\end{eqnarray*}
\begin{eqnarray*}
\nabla_{e_i}e_{m+1}&=&\nabla_{(-a_ik_ik_{i-1}\sum_{l=1}^{i-1}a_l\bar
e_l+\frac{k_i}{k_{i-1}}\bar e_i+a_ik_ik_{i-1}\bar
e_{m+1})}(\sum_{l=1}^{m}a_lk_{m}\bar e_l-k_{m}\bar e_{m+1})\\\notag
&=&-a_ik_ik_{i-1}k_{m}\sum_{l=1}^{i-1}a_l^2\nabla_{\bar e_l}\bar
e_l+\frac{a_ik_ik_{m}}{k_{i-1}}\nabla_{\bar e_i}\bar e_{i}\\\notag
&&+a_ik_ik_{i-1}k_{m}\sum_{l=1}^{i-1}a_l\nabla_{\bar e_l}\bar
e_{m+1}-\frac{k_ik_{m}}{k_{i-1}}\nabla_{\bar e_i}\bar
e_{m+1}\\\notag &=&-a_ik_ik_{i-1}k_{m}\sum_{l=1}^{i-1}a_l^2f'\bar
e_{m+1}+\frac{a_ik_ik_{m}}{k_{i-1}}f'\bar e_{m+1}\\\notag
&&-a_ik_ik_{i-1}k_{m}\sum_{l=1}^{i-1}a_lf'\bar
e_l+\frac{k_ik_{m}}{k_{i-1}}f'\bar e_i\\\notag
&=&a_ik_ik_{i-1}k_{m}(\frac{1}{k_{i-1}^2}-\sum_{l=1}^{i-1}a_l^2)f'\bar
e_{m+1},\\\notag &&-a_ik_ik_{i-1}k_{m}\sum_{l=1}^{i-1}a_lf'\bar
e_l+\frac{k_ik_{m}}{k_{i-1}}f'\bar e_i\\\notag
&=&k_{m}f'[a_ik_ik_{i-1}\bar e_{m+1}+\frac{k_i}{k_{i-1}}\bar
e_i-a_ik_ik_{i-1}\sum_{l=1}^{i-1}a_l\bar e_l]\\\notag &=&
k_{m}f'e_i,
\end{eqnarray*}
\begin{eqnarray*}
\nabla_{e_{m+1}}e_{m+1}&=&\nabla_{\sum_{l=1}^{m}a_lk_{m}\bar
e_l-k_{m}\bar e_{m+1}}(\sum_{l=1}^{m}a_lk_{m}\bar e_l-k_{m}\bar
e_{m+1})\\\notag &=&\sum_{l=1}^{m}a_l^2k_{m}^2\nabla_{\bar e_l}\bar
e_{l}-\sum_{l=1}^{m}a_lk_{m}^2\nabla_{\bar e_l}\bar e_{m+1}\\\notag
&=&(1-k_{m}^2)f'\bar e_{m+1}+\sum_{l=1}^{m}a_lk_{m}^2f'\bar e_l.
\end{eqnarray*}
We have some primary computation:
\begin{eqnarray*}\notag
\bar e_\alpha(\sigma)&=&\delta_{m+1}^{\alpha}f^{'},\; \bar e_\alpha
\bar e_\alpha(\sigma)=\delta_{m+1}^\alpha ff^{''},\; (\nabla_{\bar
e_i} \bar e_i)(\sigma)=(f^{'})^2,\\\notag (\nabla_{\bar e_{m+1}}
\bar e_{m+1})(\sigma)&=&0,\; e_i(\sigma)=
a_ik_ik_{i-1}f^{'},\;e_{m+1}(\sigma)=-k_{m}f^{'},\\\notag
 e_ie_{m+1}(\sigma)&=&-a_ik_ik_{i-1}k_{m}ff^{''},\;
 e_{m+1}e_{m+1}(\sigma)=k_{m}^2ff^{''},\\\notag
(\nabla_{e_{m+1}}e_{m+1})(\sigma)&=&(1-k_{m}^2)(f^{'})^2.
\end{eqnarray*}
A further computation gives
\begin{eqnarray*} \notag
\Delta\sigma&=&\sum_{\alpha=1}^{m+1}[ \bar e_\alpha \bar
e_\alpha(\sigma)-(\nabla_{\bar e_\alpha}\bar
e_\alpha)(\sigma)]=ff^{''}-m(f^{'})^2,\\\notag {\rm
H}ess\,(\sigma)(\xi,\xi)&=&e_{m+1}e_{m+1}(\sigma)-(\nabla_{e_{m+1}}e_{m+1})(\sigma)\\\notag
&=&k_{m}^2ff^{''}-(1-k_{m}^2)(f^{'})^2,\\\notag {\rm
grad}\,\sigma&=&f^{'}\bar e_{m+1},|{\rm
grad}\,\sigma|^2=(f^{'})^2,\;\;\;\xi(\sigma)=-k_{m}f^{'}.
\end{eqnarray*}
\begin{eqnarray}\label{RiN}
\rm Ric(\xi, \xi)&=&\Delta\sigma+(m-1)[\rm Hess(\sigma)(\xi,
\xi)-(\xi\sigma)^2+|{\rm grad}\,\sigma|^{2}]\\\notag
&=&[1+(m-1)k_m^2]ff^{''}-m(f^{'})^2.
\end{eqnarray}

Noting that $\xi=e_{m+1}$ is the unit normal vector field we can
easily compute the components of the second fundamental form as
\begin{eqnarray}\notag
h(e_i,e_i)=\langle\nabla_{e_i}e_i,e_{m+1}\rangle=-\langle
\nabla_{e_i}e_{m+1},e_i\rangle=-k_{m}f',\\\notag
h(e_i,e_j)=\langle\nabla_{e_i}e_j,e_{m+1}\rangle=-\langle\nabla_{e_i}e_{m+1},e_j\rangle=0,\;\;i\neq
j
\end{eqnarray}
from which we conclude that each of the hyperplane
$z=\sum_{i=1}^{m}a_ix_i+c$ is a totally umbilical hypersurface in
the conformally flat space and all
principal normal curvatures are equal to
\begin{equation}\label{H}
\xi(\sigma)=-k_{m}f'=H.
\end{equation}

It follows that
\begin{eqnarray}\label{RiT}
[\rm Ric(\xi)]^T&=&(m-1)[{\rm grad}_g(\xi\sigma)-\xi(\sigma){\rm
grad}_g\,\sigma+A({\rm grad}_g \,\sigma)]\\\notag &=& (m-1){\rm
grad}_g\,H,
\end{eqnarray}
and the norm of the second fundamental form is given by
\begin{equation}\label{A2}
|A|^2=\sum_{i=1}^m\langle\nabla_{e_i}\xi,\nabla_{e_i}\xi\rangle^2=mk_{m}^2(f')^2.
\end{equation}
A further computation gives
\begin{eqnarray}\notag
e_i(H)&=&(-a_ik_ik_{i-1}\sum_{l=1}^{i-1}a_l\bar
e_l+\frac{k_i}{k_{i-1}}\bar e_i+a_ik_ik_{i-1}\bar
e_{m+1})(-k_mf')\\\notag &=&a_ik_ik_{i-1}k_m\sum_{l=1}^{i-1}a_l\bar
e_l(f')-\frac{k_ik_m}{k_{i-1}}\bar e_i(f')\\\notag
&=&a_ik_ik_{i-1}k_m\sum_{l=1}^{i-1}a_l^2ff''-\frac{a_ik_ik_m}{k_{i-1}}ff''\\\notag
&=&a_ik_ik_{i-1}k_mff''(\sum_{l=1}^{i-1}a_l^2-\frac{1}{k_{i-1}^2})=-a_ik_ik_{i-1}k_mff'',
\end{eqnarray}
\begin{eqnarray}\label{gradH}
{\rm grad}_gH&=&
\sum_{i=1}^me_i(H)e_i=-\sum_{i=1}^ma_ik_ik_{i-1}k_{m}ff''e_i,
\end{eqnarray}

\begin{eqnarray}
 e_ie_i(H)&=&-e_i(a_ik_ik_{i-1}k_{m}ff'')\\\notag
&=&-(-a_ik_ik_{i-1}\sum_{l=1}^{i-1}a_l\bar
e_l+\frac{k_i}{k_{i-1}}\bar e_i)(a_ik_ik_{i-1}k_{m}ff'')\\\notag &=&
a_i^2k_i^2k_{i-1}^2k_{m}\sum_{l=1}^{i-1}a_l\bar
e_l(ff'')-a_ik_i^2k_m\bar e_i(ff'')\\\notag &=&
a_i^2k_i^2k_{i-1}^2k_{m}\sum_{l=1}^{i-1}a_l^2
f(ff'''+f'f'')-a_i^2k_i^2k_mf(ff'''+f'f'')\\\notag &=&
a_i^2k_i^2k_{i-1}^2k_{m}f(ff'''+f'f'')(\sum_{l=1}^{i-1}a_l^2-\frac{1}{k_{i-1}^2})\\\notag
&=&-a_i^2k_i^2k_{i-1}^2k_{m}(f^2f'''+ff'f''),\\\notag
(\nabla_{e_i}e_i)(H)&=&-k_i^2k_{i-1}^2k_mf'(a_i^2\sum_{l=1}^{i-1}a_l\bar
e_l-\frac{a_i}{k_{i-1}^2}\bar
e_i)(f')\\\notag&=&-a_i^2k_i^2k_{i-1}^2k_mff'f''(\sum_{l=1}^{i-1}a_l^2-\frac{1}{k_{i-1}^2})\\\notag
&=&a_i^2k_i^2k_{i-1}^2k_mff'f'',
\end{eqnarray}
\begin{eqnarray}\label{LapH}
\triangle^MH&=&\sum_{i=1}^m[e_ie_i(H)-(\nabla^M_{e_i}e_i)(H)]\;\;\;({\rm
\;and\;by\;Gauss\;formula})\\\notag
&=&\sum_{i=1}^m[e_ie_i(H)-(\nabla_{e_i}e_i)(H)+h(e_i,e_i)\,\xi\,(H)]\\\notag
&=&-\sum_{i=1}^ma_i^2k_i^2k_{i-1}^2k_m[(2-m)ff'f''+f^2f''']\\\notag&=&-(1-k_m^2)k_m[(2-m)ff'f''+f^2f'''],
\end{eqnarray}
where in obtaining the last equality we have used the identity
\begin{equation}
\sum_{i=1}^ma_i^2k_i^2k_{i-1}^2=(1-k_m^2)
\end{equation}
which can be proved by mathematical induction on $m\ge 2$. In fact,
\begin{eqnarray}\notag
a_1^2k_1^2+a_2^2k_2^2k_1^2&=&k_1^2(a_1^2+a_2^2k_2^2)=k_1^2(a_1^2+\frac{a_2^2}{1+a_1^2+a_2^2})\\\notag
&=&k_1^2\frac{(1+a_1^2)(a_1^2+a_2^2)}{1+a_1^2+a_2^2}=\frac{(a_1^2+a_2^2)}{1+a_1^2+a_2^2}=1-k_2^2,
\end{eqnarray}
so the statement is true for $m=2$. Now suppose the identity is true for $m=t\ge 2$, i.e.,
$\sum_{i=1}^ta_i^2k_i^2k_{i-1}^2=1-k_t^2$, we want to show that it is also true for $m=t+1$. Indeed, for $m=t+1$, we have
\begin{eqnarray}\notag
\sum_{i=1}^{t+1}a_i^2k_i^2k_{i-1}^2&=&\sum_{i=1}^ta_i^2k_i^2k_{i-1}^2+a_{t+1}^2k_{t+1}^2k_t^2\\\notag
&=&1-k_t^2+a_{t+1}^2k_{t+1}^2k_t^2=1+k_t^2(a_{t+1}^2k_{t+1}^2-1)\\\notag
&=&1+k_t^2(\frac{a_{t+1}^2}{1+\sum_{i=1}^{t+1}a_i^2}-1)=1-k_{t+1}^2.
\end{eqnarray}

Substituting (\ref{RiN}), (\ref{H}), (\ref{RiT}), (\ref{A2}),
(\ref{gradH}), (\ref{LapH}), and $A({\rm grad}_g\,H)=H\,{\rm
grad}_g\,H$ into biharmonic equation (\ref{eq1}) we conclude that
the isometric immersion $\varphi$ is biharmonic if and only if
\begin{equation}\label{Keq}
\begin{cases}
-(1-k_m^2)f^2f^{'''}-[3-m+(2m-3)k_m^2]ff^{'}f^{''}+m(1+k_m^2)(f')^3 =0,\\
( m-4)ff'f''a_i=0, \;\;\; i=1, 2, \ldots, m.
\end{cases}
\end{equation}
This system of ordinary differential equations can be solved by finding the solutions of the second equation of (\ref{Keq}) and then substituting these solutions into the first equation. We can
solve the second equation of (\ref{Keq}) by considering the following three
cases.  \\
{\bf Case 1.}  $f'=0$ (which implies $H=-k_mf'=0$) gives the trivial
solution. In this case, $\varphi$ is actually totally geodesic since its image is a hyperplane in a space that is homothetic to Euclidean space.\\
{\bf Case 2.}  $a_iff''=0, i=1, 2, \ldots, m$ (which, together with
(\ref{gradH}), implies ${\rm grad}_g\,H=0$). In this case the first
equation of (\ref{Keq}), with the aid of $\triangle^MH=0$ and
(\ref{LapH}), can be reduced to
\begin{equation}\label{f4}
[1+(m-1)k_m^2]ff^{''}-m(1+k_m^2)(f^{'})^2=0.
\end{equation}
If $ff''=0$, then Equation (\ref{f4}) reduces to $f^{'}=0$, which
gives the trivial solution again, i.e., the hypersurfaces are
minimal. If $ff''\ne 0$, then all $a_i=0, i=1, 2, \ldots, m$, and
hence $k_m=1$. Thus, Equation (\ref{f4}) reduces to
$ff''-2(f^{'})^2=0$ which has solutions $f(z)=\frac{1}{Az+B}$, where $A, B$ are constants.\\
{\bf Case 3.} $m=4$. In this case, the biharmonic equation
(\ref{Keq}) reduces to
\begin{eqnarray}\notag
-(1-k_4^2)f^2f'''-(5k_4^2-1)ff'f''+4(1+k_4^2)(f')^3=0,
\end{eqnarray}
from which we obtain Equation (\ref{single}).\\

Summarizing the above results we obtain the theorem.
\end{proof}

\section{The generalized Chen's conjecture on biharmonic submanifolds is false }

\indent In this section, we will show that Equation (\ref{single})
has many solutions including counter examples to the generalized
Chen's conjecture.
\begin{lemma}\label{LM1}
Let $ A>0,\; B>0,\;c $ be constants,
$\mathbb{R}^5_{+}=\{(x_1,\ldots, x_4, z)\in \r^5: z>0\}$ be the
upper-half space, and $f:\mathbb{R}^5_{+}\longrightarrow \r,\;\;
f(z)=(Az+B)^{t}$ . Then, for any $t\in (0,\;1/2)$ and any $(a_1,
a_2, a_3, a_4)\in S^3\left(\sqrt{\frac{2t}{1-2t}}\;\right)$, the
isometric immersion
\begin{equation}\label{4D}
\varphi : \mathbb{R}^4\longrightarrow
\left(\mathbb{R}^5_{+},h=f^{-2}(z)[\sum_{i=1}^4{\rm d}{x_i}^{2}+{\rm
d}{z}^{2}]\right)
\end{equation}
 with $\varphi(x_1,\ldots,
x_4)=(x_1,\ldots, x_4,\sum_{i=1}^{4}a_ix_i+c)$ is proper biharmonic
into the conformally flat space.
\end{lemma}
\begin{proof}
For $f(z)=(Az+B)^{t}$, we have
$f'=tA(Az+B)^{t-1}$,\;$f''=t(t-1)A^2(Az+B)^{t-2}$,\;$f'''=t(t-1)(t-2)A^3(Az+B)^{t-3}$.
Substituting these into Equation (\ref{single}) and using the
assumption that $A> 0$ we have
\begin{equation*}
(t-1)(t-2)\sum_{i=1}^4a_i^2+(4-\sum_{i=1}^4a_i^2)t(t-1)-4t^2(2+\sum_{i=1}^4a_i^2)=0,
\end{equation*}
which is equivalent to
\begin{equation}\label {algE}
\sum_{i=1}^4a_i^2=\frac{2t}{1-2t}.
\end{equation}

Solving the inequality $\frac{2t}{1-2t}> 0$ we conclude that for any
$t\in (0,\;1/2)$ and $(a_1, a_2, a_3, a_4)\in
S^3\left(\sqrt{\frac{2t}{1-2t}}\;\right)$ will solve the equation
(\ref{algE}) and hence the biharmonic equation (\ref{single}). Thus,
we obtain the lemma.
\end{proof}
\begin{example} Let $ A>0,\; B>0,\; c $ be constants. Let $t=\frac{1}{6}$ and $(\frac{\sqrt{2}}{4},
\frac{\sqrt{2}}{4}, \frac{\sqrt{2}}{4}, \frac{\sqrt{2}}{4})\in
S^3\left(\frac{\sqrt{2}}{2}\;\right)$, then we have a proper
biharmonic isometric immersion $\varphi:\mathbb{R}^4\longrightarrow
(\mathbb{R}^5_{+},h=(Az+B)^{-\frac{1}{3}}(\sum_{i=1}^4{\rm
d}{x_i}^{2}+{\rm d}{z}^{2}))$ with $\varphi(x_1,\ldots,
x_4)=(x_1,\ldots, x_4,\frac{\sqrt{2}}{4}(x_1+x_2+x_3+x_4)+c)$.
\end{example}

\begin{lemma}\label{LM2}
For constant $A>0,\; B>0$ and $t\in (0,1)$, the conformally flat
space $(\mathbb{R}_{+}^{5},h=(Az+B)^{-2t}(\sum_{i=1}^4{\rm
d}{x_i}^{2}+{\rm d}{z}^{2}))$ has negative sectional curvature.
\end{lemma}
\begin{proof}
Let $f(z)=(Az+B)^t$. Then, as in the proof of Theorem \ref{f1}, we
use ${\bar e}_i=f(z)\partial_i, i=1,\ldots, 5$ to denote the
orthonormal frame on
$(\mathbb{R}_{+}^{5},h=(Az+B)^{-2t}(\sum_{i=1}^4{\rm
d}{x_i}^{2}+{\rm d}{z}^{2}))$. Let $P$ be a plan section at any
point and suppose that $P$ is spanned by an orthonormal basis $X,
Y$. Then, we have $X=\sum_{i=1}^{5}a_i \bar e_i, Y=\sum_{i=1}^{5}b_i
\bar e_i$. Using sectional  curvature relation (\ref{QL}) and the
fact that the sectional curvature $\bar K(p)$ of
$(\mathbb{R}_{+}^{5},\bar h=\sum_{i=1}^4{\rm d}{x_i}^{2}+{\rm
d}{z}^{2})$ vanishes identically we find the sectional curvature of
the conformally flat space to be
\begin{eqnarray*}\notag
K(P)&=&(h(\nabla_{X}\nabla \sigma,X)+h(\nabla_{Y}\nabla
\sigma,Y))+(|\nabla\sigma|^{2}-(X\sigma )^{2}-(Y\sigma)^{2})\\\notag
&=&X(X\sigma)+Y(Y\sigma)-(\nabla_XX)(\sigma)-(\nabla_YY)(\sigma)+(|\nabla\sigma|^{2}-(X\sigma
)^{2}-(Y\sigma)^{2}),
\end{eqnarray*}
where $\sigma=\ln\,f(z)$. A straightforward computation gives

\begin{eqnarray*}\notag
X\sigma&=&\sum_{i=1}^{5}a_i \bar e_i(\sigma)=a_5f',\\\notag
X(X\sigma)&=&\sum_{i=1}^{5}a_i \bar e_i(a_5f')=\sum_{i=1}^{5}a_i
\bar e_i(a_5)f'+a_5^2ff'',
\\\notag
\nabla_XX&=&\nabla_{\sum_{i=1}^{5}a_i \bar e_i}(\sum_{j=1}^{5}a_j
\bar e_j)\\\notag &=&\sum_{i=1}^{5}a_i\bar e_i(\sum_{j=1}^{5}a_j)
\bar e_j+\sum_{i=1}^{5}\sum_{j=1}^{5}a_ia_j\nabla_{\bar e_i}\bar
e_j\\\notag &=&\sum_{i=1}^{5}a_i\bar e_i(\sum_{j=1}^{5}a_j) \bar
e_j+\sum_{i=1}^{4}a_i^2\nabla_{\bar e_i}\bar
e_i+\sum_{i=1}^{4}a_ia_5\nabla_{\bar e_i}\bar e_5\\\notag
&=&\sum_{i=1}^{5}a_i\bar e_i(\sum_{j=1}^{5}a_j) \bar
e_j+\sum_{i=1}^{4}a_i^2f'\bar e_5-\sum_{i=1}^{4}a_ia_5f'\bar
e_i\\\notag (\nabla_XX)(\sigma)&=&\sum_{i=1}^{5}a_i\bar
e_i(\sum_{j=1}^{5}a_j) \bar e_j(\sigma)+\sum_{i=1}^{4}a_i^2f'\bar
e_5(\sigma)-\sum_{i=1}^{4}a_ia_5f'\bar e_i(\sigma)\\\notag
&=&\sum_{i=1}^{5}a_i\bar e_i(a_5)f'+\sum_{i=1}^{4}a_i^2(f')^2.
\\\notag
X(X\sigma)-(\nabla_XX)(\sigma)&=&a_5^2ff''-\sum_{i=1}^{4}a_i^2(f')^2.
\end{eqnarray*}
Similarly, we have
\begin{eqnarray*}\notag
Y(Y\sigma)&=&b_5f', (\nabla_YY)(\sigma)=\sum_{i=1}^{5}b_i\bar
e_i(b_5)f'+\sum_{i=1}^{4}b_i^2(f')^2\\\notag
Y(Y\sigma)-(\nabla_YY)(\sigma)&=&b_5^2ff''-\sum_{i=1}^{4}b_i^2(f')^2,
\end{eqnarray*}
from which we have
\begin{eqnarray*}\notag
K(P)&=&X(X\sigma)+Y(Y\sigma)-(\nabla_XX)(\sigma)-(\nabla_YY)(\sigma)+(|\nabla\sigma|^{2}-(X\sigma
)^{2}-(Y\sigma)^{2})\\\notag
&=&a_5^2ff''+b_5^2ff''-\sum_{i=1}^{4}a_i^2(f')^2-\sum_{i=1}^{4}b_i^2(f')^2+(f')^2-(a_5f')^2-(b_5f')^2\\\notag
&=&(a_5^2+b_5^2)ff''+[1-(\sum_{i=1}^{4}a_i^2+a_5^2)-(\sum_{i=1}^{4}b_i^2+b_5^2)](f')^2\\\notag
&=&(a_5^2+b_5^2)ff''-(f')^2.
\end{eqnarray*}
For $f=(Az+B)^{t}$, we have
$f'=tA(Az+B)^{t-1},f''=t(t-1)A^2(Az+B)^{t-2}$, so
\begin{eqnarray*}\notag
K(P)&=&(a_5^2+b_5^2)ff''-(f')^2\\\notag &=&
A^2(Az+B)^{t-2}[(a_5^2+b_5^2)t(t-1)-t^2],
\end{eqnarray*}
which is strictly negative since $[(a_5^2+b_5^2)t(t-1)-t^2]<0$
 for $0<t<1$, and $A^2(Az+B)^{t-2}>0$ for $z>0$. From this we obtain the
 Lemma.
\end{proof}
Combining Lemma \ref{LM1} and Lemma \ref{LM2} we have
\begin{theorem}\label{MT2}
Let $ A>0,\; B>0,\;c $ be constants,
$\mathbb{R}^5_{+}=\{(x_1,\ldots, x_4, z)\in \r^5: z>0\}$ be the
upper-half space, and $f:\mathbb{R}^5_{+}\longrightarrow \r,\;\;
f(z)=(Az+B)^{t}$ . Then, for any $t\in (0,\;1/2)$ and any $(a_1,
a_2, a_3, a_4)\in S^3\left(\sqrt{\frac{2t}{1-2t}}\;\right)$, the
isometric immersion
\begin{equation}\label{5D}
\varphi : \mathbb{R}^4\longrightarrow
\left(\mathbb{R}^5_{+},h=f^{-2}(z)[\sum_{i=1}^4{\rm d}{x_i}^{2}+{\rm
d}{z}^{2}]\right)
\end{equation}
 with $\varphi(x_1,\ldots,
x_4)=(x_1,\ldots, x_4,\sum_{i=1}^{4}a_ix_i+c)$ gives a proper
biharmonic hypersurface into the conformally flat space with
strictly negative sectional curvature. These provide infinitely many
counter examples to the generalized Chen's conjecture on biharmonic
submanifolds.
\end{theorem}
The following corollary can be used to construct proper biharmonic
submanifolds of any codimension in a nonpositively curved manifold.
\begin{corollary}
For any positive integer $k$, there exists a proper biharmonic
submanifold of codimension $k$ in a nonpositively curved space.
Thus, the generalized Chen's conjecture is false.
\end{corollary}
\begin{proof}
Let \begin{equation}\label{4D} \varphi : \mathbb{R}^4\longrightarrow
\left(\mathbb{R}^5_{+},h=f^{-2}(z)[\sum_{i=1}^4{\rm d}{x_i}^{2}+{\rm
d}{z}^{2}]\right)
\end{equation}
 with $\varphi(x_1,\ldots,
x_4)=(x_1,\ldots, x_4,\sum_{i=1}^{4}a_ix_i+c)$ be one of the proper
biharmonic hypersurface given in Theorem \ref{MT2} and $\psi:
\mathbb{R}^n\longrightarrow \mathbb{R}^n\times
\mathbb{R}^{k-1}\equiv (\mathbb{R}^{n+k-1}, h_0)$ with
$\psi(y)=(y,0)$ be the totally geodesic embedding of a subspace into
a Euclidean space. Then, the isometric embedding $\phi:
\mathbb{R}^4\times \mathbb{R}^n\longrightarrow
\left(\mathbb{R}^5_{+}\times \mathbb{R}^{n+k-1}, h+h_0\right)$ with
$\phi(x,y)=(\varphi(x), \psi(y))$ gives a submanifold of codimension
$k$. Since $\phi$ is biharmonic with respect to each variable
separately and it is proper biharmonic with respect to $x$-variable
by Theorem \ref{MT2}, we can use Proposition 2.1 in \cite{Ou2} to
conclude that $\phi$ is a proper biharmonic embedding. Thus, the
image of $\phi$ provides a proper biharmonic submanifold of
codimension $k$. Since, by Lemma \ref{LM2}, the conformally flat
space $(\mathbb{R}_{+}^{5},h=(Az+B)^{-2t}(\sum_{i=1}^4{\rm
d}{x_i}^{2}+{\rm d}{z}^{2}))$ has negative sectional curvature and
the Euclidean space $(\mathbb{R}^{n+k-1}, h_0)$ has zero curvature,
their product $\left(\mathbb{R}^5_{+}\times \mathbb{R}^{n+k-1},
h+h_0\right)$ gives a space of nonpositive curvature. Thus, we
obtain the corollary.
\end{proof}
\begin{ack}
Both authors are grateful to Cezar Oniciuc for a careful reading and
checking the computations in the paper. We especially want to thank
him for pointing out a mistake in the computation of the Laplacian
of the mean curvature of the hypersurface in the previous version.
\end{ack}

\end{document}